\documentclass[twoside,leqno,10pt]{article}

\usepackage{graphics,ifthen}
\usepackage{latexsym}
\usepackage{mathrsfs}
\usepackage{amsmath}
\usepackage{amssymb}
\newcommand{\Pcr}{\mathscr{P}}

\newcommand{\Scr}{\mathscr{S}}

\newcommand{\Xcr}{\mathscr{X}}

\newcommand{\E}{\mathbb{E}}

\newcommand{\X}{\mathbb{X}}
\newcommand{\Pe}{\mathbb{P}}

\begin{document}
\input{latexP.sty}
\input{referencesP.sty}
\input epsf.sty

\def\ind{\stackrel{\mathrm{ind}}{\sim}}
\def\iid{\stackrel{\mathrm{iid}}{\sim}}
\def\Prodi{\mathop{{\lower9pt\hbox{\epsfxsize=15pt\epsfbox{pi.ps}}}}}
\def\prodi{\mathop{{\lower3pt\hbox{\epsfxsize=7pt\epsfbox{pi.ps}}}}}

\def\Definition{\stepcounter{definitionN}\
    \Demo{Definition\hskip\smallindent\thedefinitionN}}
\def\EndDefinition{\EndDemo}
\def\Example#1{\Demo{Example [{\rm #1}]}}
\def\EndExample{\qed\EndDemo}
\def\Category#1{\centerline{\Heading #1}\rm}
\
\def\e{\text{\hskip1.5pt e}}
\newcommand{\eps}{\epsilon}
\newcommand{\proof}{\noindent {\bf Proof:\ }}
\newcommand{\remarks}{\noindent {\bf Remarks:\ }}
\newcommand{\note}{\noindent {\bf Note:\ }}
\newcommand{\examp}{\noindent {\bf Example:\ }}
\newcommand{\Lower}[2]{\smash{\lower #1 \hbox{#2}}}
\newcommand{\ben}{\begin{enumerate}}
\newcommand{\een}{\end{enumerate}}
\newcommand{\bi}{\begin{itemize}}
\newcommand{\ei}{\end{itemize}}
\newcommand{\hp}{\hspace{.2in}}

\newtheorem{lw}{Proposition 3.1, Lo and Weng (1989)}
\newtheorem{thm}{Theorem}[section]
\newtheorem{defin}{Definition}[section]
\newtheorem{prop}{Proposition}[section]
\newtheorem{lem}{Lemma}[section]
\newtheorem{cor}{Corollary}[section]
\newcommand{\rb}[1]{\raisebox{1.5ex}[0pt]{#1}}
\newcommand{\mc}{\multicolumn}
\def\Beta{\text{Beta}}
\def\Dir{\text{Dirichlet}}
\def\DP{\text{DP}}
\def\P{{\bf p}}
\def\fhat{\widehat{f}}
\def\GA{\text{gamma}}
\def\ind{\stackrel{\mathrm{ind}}{\sim}}
\def\iid{\stackrel{\mathrm{iid}}{\sim}}
\def\K{{\bf K}}
\def\min{\text{min}}
\def\N{\text{N}}
\def\p{{\bf p}}
\def\U{{\bf U}}
\def\T{{\bf T}}
\def\m{{\bf m}}
\def\X{{\bf X}}
\def\Y{{\bf Y}}
\def\tps{\mbox{\scriptsize ${\theta H}$}}   
\def\ups{\mbox{\scriptsize ${P_{\theta}(g)}$}}   
\def\vps{\mbox{\scriptsize ${\theta}$}}   
\def\vups{\mbox{\scriptsize ${\theta >0}$}}   
\def\hps{\mbox{\scriptsize ${H}$}}   
\def\rps{\mbox{\scriptsize ${(\theta+1/2,\theta+1/2)}$}}   
\def\sps{\mbox{\scriptsize ${(1/2,1/2)}$}}   

\newcommand{\reals}{{\rm I\!R}}
\newcommand{\PR}{{\rm I\!P}}
\def\Z{{\bf Z}}
\def\yy{{\mathcal Y}}
\def\rr{{\mathcal R}}
\def\BP{\text{beta}}
\def\ts{\tilde{t}}
\def\js{\tilde{J}}
\def\gs{\tilde{g}}
\def\fs{\tilde{f}}
\def\ys{\tilde{Y}}
\def\ps{\tilde{\mathcal {P}}}

\def\Report{Lancelot F. James}
\def\Author{NTR mixture models}
\pagestyle{myheadings}
\markboth{\Author}{\Report}
\thispagestyle{empty}

\bct\Heading Spatial Neutral to the Right Species Sampling Mixture
Models\lbk\lbk\smc Lancelot F. James\footnote{ \eightit AMS 2000
subject classifications.
               \rm Primary 62G05; secondary 62F15.\\
\eightit Corresponding authors address.
                \rm The Hong Kong University of Science and Technology,
Department of Information Systems and Management,
Clear Water Bay, Kowloon,
Hong Kong.
\rm lancelot\at ust.hk\\
\indent\eightit Keywords and phrases.
                \rm
          Chinese Restaurant process,
          Dirichlet process,
          L\'evy processes,
          Neutral to the right processes,
          Species sampling models.
          }
\lbk\lbk \BigSlant The Hong Kong University of Science and
Technology\rm \lbk
\ect \Quote  This paper describes briefly how one may utilize a
class of species sampling mixture models derived from
Doksum's~(1974) neutral to the right processes. For practical
implementation we describe an ordered/ranked variant of the
generalized weighted Chinese restaurant process. \EndQuote
\rm

\section{Introduction}
The field of Bayesian nonparametric statistics essentially
involves the idea of assigning prior and posterior distributions
over spaces of probability measures or more general measures. That
is, similar to the classical parametric Bayesian idea of assigning
priors to an unknown parameter, say $\theta,$ which lies in a
Euclidean space, one views, for instance, an unknown cumulative
distribution function, say $F(t)$, as being a stochastic process.
More generally for an unknown probability measure $P$, a Bayesian
views it as a random probability measure. This is currently a
well-developed and active area of research that has links to a
variety of areas where L\'evy and more general random processes
are commonly used. However, as discussed in~Doksum and
James~(2004), in the late 1960's, noting the high activity and
advance in nonparametric statistics, David Blackwell and others
wondered how one could assign priors which were both flexible and
tractable. Arising from these questions were two viable answers
which till this day remain at the cornerstone of Bayesian
nonparametric statistics.

Ferguson~(1973, 1974) proposed the use of a Dirichlet process
prior[see also Freedman~(1963)]. For this prior if $P$ is a
probability on some space $\Xcr,$ and $(B_{1},\ldots, B_{k})$ is a
measurable partition of $\Xcr,$ then $P(B_{1}),\ldots, P(B_{k})$
has a Dirichlet distribution. Moreover, the posterior distribution
of $P$ given a sample $\X=(X_{1},\ldots, X_{n})$ is also a
Dirichlet process. For a specified probability measure $H$ and a
scalar $\theta>0$, one can say that $P:\overset {d}=P_{\theta H}$
is a Dirichlet process with shape parameter $\theta H$, if the
Dirichlet distributions discussed above have parameters given by
$\E[P(A_{i})]=\theta H(A_{i}).$ Following this, Doksum (1974)
introduced the class of Neutral to the Right (NTR) random
probability measures on the real line. For these models if $P$ is
a distribution on the real line, then for each partition
$B_{1},\ldots, B_{k}$, with $B_{j}=(s_{j-1},s_{j}],$ $j=1,\ldots,
k$, $s_{0}=-\infty,s_{k}=\infty$, $s_{i}<s_{j}$ for $i<j$;
$P(B_{1}),\ldots, P(B_{k})$ is such that $P(B_{i})$ has the same
distribution as $V_{i}\prod_{j=1}^{i-1}(1-V_{j})$, where
$V_{1},\ldots, V_{2},\ldots$ is a collection of independent
non-negative random variables. This represents a remarkably rich
choice of models defined by specifying different distributions for
the $V_{i}.$ Notably if $V_{i}$ is chosen to be beta random
variable with parameters $(\alpha_{i},\beta_{i})$ and
$\beta_{i}=\sum_{j=1}^{k-1}\alpha_{j}$, then this gives the
Dirichlet process as described in Doksum~(1974). Doksum~(1974)
shows that if $P$ is a NTR distribution then the posterior
distribution of $P$ give a sample $X_{1},\ldots, X_{n}$ is also an
NTR. Subsequently, Ferguson and Phadia~(1979), showed that this
type of conjugacy property extends to the case of right censored
survival models. This last fact coupled with the subsequent
related works of Hjort~(1990), Kim~(1999), Lo~(1993) and Walker
and Muliere~(1997) have popularized the usage of NTR processes in
models related to survival and event history analysis.

Despite these attractive points, the usage of NTR processes in
more complex statistical models, such as mixture models, has been
notably absent. This is in contrast to the Dirichlet process
which, coupled with the advances in MCMC and other computational
procedures, is regularly used in nonparametric or semi-parametric
statistical models. The theoretical framework for Dirichlet
process mixture models can be traced back to Lo~(1984) who
proposed to model a density as a  convolution mixture model of a
known kernel density $K(y|x)$ and a Dirichlet process $P$ as, \Eq
f(y|P)=\int_{\Xcr}K(y|x)P(dx)\label{LoDen}. \EndEq This may be
equivalently expressed in terms of a missing data model where for
a sample $\Y=(Y_{1},\ldots, Y_{n})$ based on~\mref{LoDen}, one has
$Y_{1},\ldots, Y_{n}|\X,P$ are such that $Y_{i}$ are independent
with distributions $K(\cdot|X_{i})$, $X_{i}|P$ are iid $P$ and $P$
is a Dirichlet process. It is clear that the description of the
posterior distribution of $P$ and related quantities is much more
complex than in the setting discussed in Ferguson~(1973). However,
Lo~(1984) shows that its description is facilitated by the
descriptions of the posterior distribution of $P|\X$, given by
Ferguson~(1973) and the exchangeable marginal distribution of $\X$
discussed in Blackwell and MacQueen~(1973). Blackwell and Macqueen
describe the distribution via what is known as the
Blackwell-MacQueen P\'olya urn scheme where $\Pe(X_{1}\in A)=H(A)$
and for $n>1$ \Eq \Pe(X_{n}\in \cdot|X_{1},\ldots,
X_{n-1})=\frac{\theta}{\theta+n-1}H(\cdot)+\frac{1}{\theta+n-1}\sum_{j=1}^{n-1}\delta_{X_{i}}(\cdot)
\label{BMQ}.\EndEq

Note that~\mref{BMQ} clearly indicates that there are ties among
$(X_{1},\ldots, X_{n})$ and that the $n(\p)\leq n$ unique values,
say $X^{*}_{1},\ldots,X^{*}_{n(\p)}$ are iid with common
distribution $H.$ Letting $\p=\{C_{1},\ldots, C_{n(\p)}\}$ denote
a partition of the integers $\{1,\ldots, n\}$, where one can write
$C_{j}=\{i:X_{i}=X^{*}_{j}\}$, with size $n_{j}=|C_{j}|$ for
$j=1,\ldots, n(\p).$ This leads to the following important
description of the distribution of $\X,$
$$
\pi(d\X|\theta
H)=\textsc{PD}(\p|\theta)\prod_{j=1}^{n(\p)}H(dX^{*}_{j})
$$
where
$$PD(\p|\theta)=\frac{\theta^{n(\p)}\Gamma(\theta)}{\Gamma(\theta+n)}\prod_{j=1}^{n(\p)}(n_{j}-1)!:=p_{\theta}(n_{1},\ldots,n_{n(\p)})
$$ is a variant of Ewens sampling formula[see Ewens~(1972) and
Antoniak~(1974)], often called the Chinese restaurant process. It
can be interpreted as
$\Pe(C_{1},\ldots,C_{n(\p)})=p_{\theta}(n_{1},\ldots,n_{n(\p)})$
where $p_{\theta}$, being symmetric in its arguments, is the most
notable example of an  \emph{exchangeable partition probability
function}(EPPF)~[see Pitman~(1996)]. It is easily seen that a
Dirichlet Process with shape $\theta H$ is characterized by the
pair $(p_{\theta},H).$ Letting $p(n_{1},\ldots, n_{k})$, for
$n(\p)=k,$ denote an arbitrary EPPF, Pitman~(1996) shows that the
class of random probability measures whose distribution is
completely determined by the pair $(p,H)$ must correspond to the
class of \emph{species sampling random probability measures}.
General species sampling random probability measures constitute
all random probability measures that can be represented as \Eq
P(\cdot)=\sum_{i=1}^{\infty}P_{i}\delta_{Z_{i}}(\cdot)+(1-\sum_{k=1}^{\infty}P_{k})H(\cdot)
\label{SSM}\EndEq where $0\leq P_{i}<1$ are random weights such
that $0<\sum_{i=1}^{\infty}P_{i}\leq 1$, independent of the
$Z_{i}$ which are iid with some non-atomic distribution $H$.
Furthermore the law of the $(P_{i})$ is determined by the EPPF
$p.$ Noting these points Ishwaran and James~(2003) described the
class of species sampling mixture models by replacing a Dirichlet
process in~\mref{LoDen} with $P$ specified by~\mref{SSM}. See also
M\"uller and Quintana~(2004).

Except for the special case of the Dirichlet process, NTR
processes are not species sampling models and this is one of the
factors which makes analysis a bit more difficult. Nonetheless,
  James~(2003, 2006) was able to extend the definition of NTR processes
  to a class of random probability measures on more general
  spaces, which he called Spatial NTR processes. Additionally a
  tractable description of the marginal distribution of this class
  of models was obtained. These two ingredients then allow for the
  implementation of NTR mixture models. Our goal in this note is
  not to describe the mechanisms for a full-blown NTR mixture
  model, as this requires much more overhead, but
  rather mixture models based on species sampling models which are
  derived from NTR processes. James~(2003, 2006) introduced and calls these \emph{NTR species
  sampling models}.  Quite specifically, though the NTR processes
  are not species sampling models they produce EPPF's $p$ that, along with the specification of $H$, are uniquely associated with
an
 NTR species sampling model. This produces a very rich and flexible class of
random
  priors that are a bit simpler analytically than NTR processes. An interesting fact is that this class contains the two-parameter $(\alpha,\theta)$ Poisson-Dirichlet random probability
  measures for parameters $0\leq \alpha<1$ and $\theta>0.$ That is
  the Dirichlet process and a class of random probabilities defined by normalizing a stable
  law process and further power tempering the stable law
  distribution, which are discussed in Pitman~(1996) and Pitman and Yor~(1997). Implementations of these latter models, being quite special, may be treated by computational procedures involving random partitions discussed in Ishwaran and James~(2003) or by the methods in Ishwaran and James~(2001).
  Here we will
  discuss a ranked weighted Chinese restaurant procedure which applies more
  generally.
\Section{NTR and related processes}Let $F(t)$ denote an NTR
cumulative distribution function on the positive real line.
Additionally, let $S(t)=1-F(t)$ denote a survival function.
 Doksum~(1974) Theorem 3.1 shows that $F$ is an NTR
process if and only if it can be represented as \begin{equation}
F(t)=1-{\mbox e}^{-Y(t)} \label{expref}
\end{equation}
where $Y(t)$ is an independent increment process which is
non-decreasing and right continuous almost surely and furthermore
$\lim_{t\rightarrow \infty}Y(t)=\infty$ and $\lim_{t\rightarrow
-\infty}Y(t)=0$ almost surely. In other words $Y$ belongs to the
class of positive L\'evy processes. We shall suppose hereafter
that $T$ is a positive random variable such that given $F$ its
distribution function is $F$ where $F$ is an NTR process. Then $T$
has an interpretation as a survival time with ``conditional"
survival distribution $S(t)=1-F(t):=P(T>t|F)$. It is evident from
\mref{expref} that the distribution of $F$ is completely
determined by the law of $Y$ which is determined by its Laplace
transform
$$
E\[{\mbox e}^{-\omega Y(t)}\]={\mbox
e}^{-\int_{0}^{t}\phi(\omega|s)\Lambda_{0}(ds)}:=E[(S(t))^{\omega}]
$$
where $\phi(\omega|s)$ is equal to \Eq \int_{0}^{\infty}(1-{\mbox
e}^{-v\omega})\tau(dv|s)=\int_{0}^{1}(1-{(1-u)}^{\omega})\rho(du|s)=\int_{0}^{1}\omega{(1-u)}^{\omega-1}
\[\int_{u}^{1}\rho(dv|s)\]du\label{room}\EndEq
$\tau$ and $\rho$ are L\'evy densities on $[0,\infty]$ and $[0,1]$
respectively which are in correspondence via the mapping
$y\rightarrow 1-e^{-y}$. Without loss of generality we shall
assume that $\int_{0}^{1}u\rho(du|s)=1$ for each fixed $s$, which
implies that $\phi(\omega|s)=1$. Hence we have that
$$
E[S(t)]={\mbox e}^{-\Lambda_{0}(t)}=1-F_{0}(t)
$$
where $F_{0}$ represents one's prior belief about the true
distribution and $\Lambda_{0}(dt)=F_{0}(dt)/S_{0}(t-)$ is its
corresponding cumulative hazard with
$S_{0}(t-)=1-F_{0}(t-)=\Pe(T\ge t).$

Note that for each fixed $s$, $\phi(\omega|s)$ corresponds to the
log Laplace transform of an infinitely-divisible random variable.
It follows that different specifications for $\tau$ or
equivalently $\rho$ lead to different NTR processes. When $\tau$
and $\rho$ do not depend on $s$, $F$, $Y$ and all relevant
functionals are said to be \emph{homogeneous}. We also apply this
name to $\tau$ and $\rho$. Additionally $\phi(\omega|s)$
specializes to $$ \phi(\omega):=\int_{0}^{\infty}(1-{\mbox
e}^{-v\omega})\tau(dv)=\int_{0}^{1}(1-{(1-u)}^{\omega})\rho(du)=\int_{0}^{1}\omega{(1-u)}^{\omega-1}
\[\int_{u}^{1}\rho(dv)\]du.$$
Consider now the cumulative hazard process of $F$, say $\Lambda,$
defined by $\Lambda(dt)=F(dt)/S(t-).$ The idea of Hjort~(1990) was
to work directly with  $\Lambda$ rather than $F$. He showed
importantly that if one specified $\Lambda$ to be a positive
completely random measure on $[0,1]$, whose law is specified by
the Laplace transform $$\E[{\mbox e}^{-\omega\Lambda(t)}]={\mbox
e}^{-\int_{0}^{t}\psi(\omega|s)\Lambda_{0}(ds)}$$ where
$\psi(\omega|s):=\int_{0}^{1}(1-{\mbox e}^{-u\omega})\rho(du|s),$
then $F$ and $S$ must be NTR processes specified by~\mref{room}.
James~(2003, 2006) shows that one can extend the definition of an
NTR process to a spatial NTR process on $[0,\infty]\times \Xcr$ by
working with the concept of a random hazard measure, say
$\Lambda_{H}(dt,dx).$ $\Lambda_{H}$ is a natural extension of
$\Lambda$ in the sense that $\Lambda_{H}(dt,\Xcr)=\Lambda(dt)$ and
is otherwise specified by replacing the intensity
$\rho(du|s)\Lambda_{0}(ds)$ by $\rho(du|s)\Lambda_{0}(ds,dx)$,
where,
$$
\Lambda_{0}(ds,dx)=H(dx|s)\Lambda_{0}(ds)
$$
is a hazard measure and $H(\cdot|s)$ may be interpreted as the
conditional distribution of $X|T=s.$ A Spatial NTR process
(SPNTR)~is then defined as \Eq
P_{S}(dt,dx)=S(t-)\Lambda_{H}(dt,dx)\label{SPNTR}\EndEq The SPNTR
in~\mref{SPNTR} has marginals such that $P_{S}(dt,d\Xcr)=F(dt)$ is
an NTR and \Eq
P_{S}([0,\infty),dx)=\int_{0}^{\infty}S(t-)\Lambda_{H}(ds,dx),\label{MSPNTR}\EndEq
represents an entirely new class of random probability measures.

\subsection{NTR species sampling models}
NTR species sampling models arise as a special case
of~\mref{MSPNTR} by setting $H(dx|s):=H(dx).$ Here we will further
work only with the class of homogeneous processes and hence we
will additionally choose $\rho(du|s)=\rho(du).$ Thus an NTR
species sampling model is of the form
$$
P_{\rho,H}(dx)=\int_{0}^{\infty}S(s-)\Lambda_{H}(ds,dx)=\sum_{k=1}^{\infty}P_{k}\delta_{Z_{k}}(dx).
$$
Furthermore, if $P\overset{d}=P_{\rho,H}$ the we denote its law as
$\Pcr(\cdot|\rho,H).$ It follows that for practical usage in
mixture models one needs a tractable description of the
corresponding EPPF, say $p_{\rho}.$ However, before we do that we
will need to introduce additional notation which connects
$p_{\rho}$ with the NTR process. If we suppose that $X_{1},\ldots,
X_{n}|P_{\rho,H}$ are iid with distribution $P_{\rho,H}$, then
these points come from a description of the $n$ conditionally
independent pairs $(T_{1},X_{1}),\ldots, (T_{n},X_{n})|P_{S}$
where $(T_{i},X_{i})$ are iid $P_{S}$, such that $T_{i}$ are iid
$F$, where $F$ is an NTR,  and $X_{i}$ are iid $P_{\rho,H}.$ Here
$P_{S}$ must be specified by the intensity
$\rho(du)\Lambda_{0}(ds)H(dx).$ Now if one denotes the $n(\p)$
unique pairs as $(T^{*}_{j},X^{*}_{j})$ for j=1,\ldots, n(\p),
then one may simply set each $C_{j}=\{i:T_{i}=T^{*}_{j}\}.$
Furthermore we define $T_{({1:n})}>T_{({2:n})}>\ldots
>T_{({n(\p):n})}>0$ to be the ordered values of the unique values
$(T^{*}_{j})_{j\leq n(\p)}.$ Hence we can define $\p$ by setting
$C_{j}:=\{i:T_{i}=T^{*}_{j}\}$, and define
$\m=\{D_{1},\ldots,D_{n(\p)}\}$ with cells
$D_{j}=\{i:T_{i}=T_{(j:n)}\}$ with cardinality $d_{j}=|D_{j}|.$ It
is evident that given a partition $\p=\{C_{1},\ldots,C_{n(\p)}\}$,
$\m$ takes its values over the symmetric group, say
$\Scr_{n(\p)}$, of all $n(\p)!$ permutations of $\p.$ Let
$R_{j-1}=\bigcup_{k=1}^{j-1}D_{k}:=\{i:T_{i}>T_{(j:n)}\}$ with
cardinality $r_{j-1}=\sum_{k=1}^{j-1}d_{k}$. Then, in terms of
survival analysis, the quantities $d_{j}$ and
$r_{j}=d_{j}+r_{j-1}$ have the interpretation as the number of
deaths at time $T_{(j:n)}$, and the number at risk at time
$T_{(j:n)},$ respectively. See James~(2006) for some further
elaboration. Now from James~(2003, 2006) it follows that \Eq
\pi_{\rho}(\p)=p_{\rho}(n_{1},\ldots, n(\p))=\sum_{\m\in
\Scr_{{n(\bf p)}}}\frac{ \prod_{j=1}^{n({\bf p})}\kappa_{d_{j},
r_{j-1}}(\rho)}{\prod_{j=1}^{n({\bf p})}
\phi(r_{j})}\label{SEPPF}\EndEq where,
$$
\kappa_{d_{j},r_{j-1}}(\rho)=\int_{0}^{1}u^{d_{j}}{(1-u)}^{r_{j-1}}\rho(du).$$
The form of the EPPF is in general not directly tractable. However
by augmentation one sees that the distribution of $\m$ is given by
\Eq \pi_{\rho}(\m)=\frac{ \prod_{j=1}^{n({\bf p})}\kappa_{d_{j},
r_{j-1}}(\rho)}{\prod_{j=1}^{n({\bf p})} \phi(r_{j})}
\label{Comp}\EndEq and has a nice product form. This suggests that
one can work with a joint distribution of $(\X,\m)$ given by
$$
\pi_{\rho}(\m)\prod_{j=1}^{n(\p)}H(dX^{*}_{j}).
$$
Related to this, James~(2006) shows that a prediction rule of
$X_{n+1}|\X,\m$ is given by $$\Pe(X_{n+1}\in
dx|\X,\m)=(1-\sum_{j=1}^{n(\p)}p_{j:n})P_{0}(dx)+
\sum_{j=1}^{n(\p)}p_{j:n}\delta_{X^{*}_{j}}(dx), $$ with
$(1-\sum_{j=1}^{n(\p)}p_{j:n})=\sum_{j=1}^{n(\p)+1}q_{j:n}$, and
where $$p_{j:n}
=\frac{\kappa_{d_{j}+1,r_{j-1}}(\rho)\prod_{l=j+1}^{n(\p)}
\kappa_{d_{l},r_{l-1}+1}(\rho)}
{\kappa_{d_{j},r_{j-1}}(\rho)\prod_{l=j+1}^{n(\p)}\kappa_{d_{l},r_{l-1}}(\rho)}
\prod_{l=j}^{n(\p)}\frac{\phi(r_{l})}{\phi(r_{l}+1)}. $$ and
$$q_{j:n}=
\frac{\kappa_{1,r_{j-1}}(\rho)}{\phi(r_{j-1}+1)}\frac{\prod_{l=j}^{n(\p)}
\kappa_{d_{l},r_{l-1}+1}(\rho)}
{\prod_{l=j}^{n(\p)}\kappa_{d_{l},r_{l-1}}(\rho)}
\prod_{l=j}^{n(\p)}\frac{\phi(r_{l})}{\phi(r_{l}+1)},$$ with
$q_{n(\p)+1:n}=\kappa_{1,n}(\rho)/\phi(n+1)$, are transition
probabilities derived from $\pi_{\rho}(\m).$ Note that in the
calculation of $\kappa_{1,r_{j-1}}(\rho)$, $r_{j-1}+1$ is to be
used rather than $r_{j}=r_{j-1}+m_{j}$. As an example, consider
the choice of a homogeneous beta process~[Hjort~(1990), see also
Ferguson~(1974), Ferguson and Phadia~(1979) and Gnedin~(2004)]
defined by
$$
\rho(du)=\theta u^{-1}(1-u)^{\theta -1}
$$
then it is easily seen that
$\phi(r_{j})=\sum_{l=1}^{r_{j}}\t/(\t+l-1)$, and it follows that
in this case $$p_{j:n}=\frac{d_{j}}{n+\t}
\prod_{l=j}^{n(\p)}\frac{\phi(r_{l})}{\phi(r_{l}+1)} {\mbox { and
  }}q_{j:n}=\frac{1}{n+\t}
\frac{1}{\sum_{i=1}^{r_{j-1}+1}1/(\t+i-1)}
\prod_{l=j}^{n(\p)}\frac{\phi(r_{l})}{\phi(r_{l}+1)}.$$

\Remark Gnedin and Pitman~(2005a) also obtained the
expressions~\mref{SEPPF} and~\mref{Comp} independent of
James~(2003, 2006), and in a different context. See James~(2006)
for more details.\EndRemark \Remark Related to this, Gnedin and
Pitman~(2005a)~[see additionally Gnedin and Pitman~(2005b)] showed
that the EPPF in~\mref{SEPPF} corresponds to that of the
two-parameter $(\alpha,\theta)$ Poisson-Dirichlet process with
parameters $0\leq \alpha<1$ and $\theta>0$ if
$\rho:=\rho_{\alpha,\theta}$ is chosen such that, $$
\int_{u}^{1}\rho_{\alpha,\theta}(dv)=\frac{\Gamma(\theta+2-\alpha)}{\Gamma(1-\alpha)\Gamma(1+\theta)}u^{-\alpha}{(1-u)}^{\theta}
.$$From this, James~(2006) deduced that $
P_{\rho_{\alpha,\theta},H}=\sum_{k=1}^{\infty}W_{k}\prod_{i=1}^{k-1}(1-W_{i})\delta_{Z_{k}}
$ where $(W_{k})$ are independent beta $(1-\alpha,
\theta+k\alpha)$ random variables independent of the $(Z_{k})$
which are iid $H$. That is a two-parameter $(\alpha,\theta)$
Poisson-Dirichlet process, for $0\leq \alpha<1$ and $\theta>0$ can
be represented as the marginal probability measure of a spatial
NTR process, as described above. See Pitman and Yor~(1997) and
Ishwaran and James~(2001) for more on the stick-breaking
representation of the two parameter Poisson-Dirichlet
process.\EndRemark \Section{NTR species sampling mixture models}
Now setting $P=P_{\rho,H}$ in~\mref{LoDen} yields a special case
of the species sampling models described in Ishwaran and
James~(2003). That is \Eq
\int_{\Xcr}K(y|x)P_{\rho,H}(dx)=\int_{\Xcr}\int_{0}^{\infty}K(y|x)S(s-)\Lambda_{H}(ds,dx)
\label{Sden}\EndEq is called an NTR species sampling models. We
look at the situation where $Y_{1},\ldots,Y_{n}|P_{\rho,H}$ are
iid with density or pmf~\mref{Sden}. This translates into the
hierarchical model,

\begin{eqnarray}
  Y_{i}|X_{i},P &\overset{ind}\sim& K(Y_{i}|X_{i}) {\mbox { for }} i=1,\ldots,n\nonumber\\
  X_{i}|P &\overset{iid}\sim& P \label{hier}\\
  P &\sim& \Pcr(\cdot|\rho,H)\nonumber
\end{eqnarray}

In principle, since we have a description of the EPPF, the
theoretical results and computational procedures described in
Ishwaran and James~(2003) apply. However as we have noted in
general  $\pi_{\rho}(\p)$ is not as simple to work with as
$\pi_{\rho}(\m)$. So here we develop results that allows us to
sample from a posterior distribution of $\m$ rather than
partitions. We summarize these results in the next proposition
\begin{prop} Suppose that one has the model specified
in~\mref{hier}. Then the following results holds \Enumerate
\item[(i)]The distribution of $X_{1},\ldots,X_{n}|\Y,\m$ is such
that the unique values $X^{*}_{j}$ for $j=1,\ldots,n(\p)$ are
conditionally independent with distributions
$$
\pi(dX^*_{j}|D_{j})\propto H(dX^*_{j})\prod_{i\in
D_{j}}K(Y_{i}|X^*_{j}).$$
\item[(ii)]The posterior distribution of $\m|\Y$ is,
$$
\pi_{\rho}(\m|\Y)\propto
\pi_{\rho}(\m)\prod_{j=1}^{n(\p)}\int_{\Xcr}\prod_{i\in
D_{j}}K(Y_{i}|x)H(dx).
$$
\item[(iii)] The posterior distribution of $\p|\Y$ is
$$
\sum_{\m\in
\Scr_{n(\p)}}\pi_{\rho}(\m|\Y)=\pi_{\rho}(\p)\prod_{j=1}^{n(\p)}\int_{\Xcr}\prod_{i\in
C_{j}}K(Y_{i}|x)H(dx).
$$\qed
 \EndEnumerate \end{prop}

From this result one can compute a Bayesian predictive density of
$Y_{n+1}|\m,\Y$ as,
$$
l(n)=f(Y_{n+1}|\m,\Y)=\[\sum_{j=1}^{n(\p)+1}q_{j:n}\]\int_{\Xcr}K(Y_{n+1}|x)H(dx)+\sum_{j=1}^{n(\p)}p_{j:n}\int_{\Xcr}K(Y_{n+1}|x)
\pi(dx|D_{j}).
$$
A Bayesian density estimate analogous to Lo (1984) is then to sum
this expression relative to the distribution of $\m|\Y$.
\begin{cor}
Consider the model in Proposition 3.1, then a Bayesian predictive
density estimator of $Y_{n+1}|\Y$ is given by
$$
\E[f(Y_{n+1}|P)|\Y] =\sum_{\p}\sum_{\m\in
S_{n(\p)}}f(Y_{n+1}|\m,\Y)\pi_{\rho}(\m|\Y)
$$\qed
\end{cor}

\subsubsection{Ordered/Ranked generalized weighted Chinese restaurant processes}
The significance of the expression for the predictive density, is
that we can use $l(n)$ in precisely the same manner as the
predictive densities given $\p,\Y$, used in Ishwaran and
James~(2003)~[see also Lo, Brunner and Chan~(1996)] to construct
computational procedures for approximating posterior quantities.
In fact, all the major computational procedures for Dirichlet
process mixture models, see for instance Escobar (1994) and
Escobar and West~(1995), utilize some type of predictive density.
Here, in analogy to the gWCR algorithms in Lo, Brunner and
Chan~(1996) and Ishwaran and James~(2003), we define a weighted
version of the {\it Ordered/Ranked generalized Chinese restaurant
process} developed in James~(2003, 2006), to approximate a draw
from $\pi_{\rho}(\m|\Y)$ as follows. For each $n\ge 1$, let
$\{D_{1:n},\ldots, D_{n(\p):n}\}$, denote a seating configuration
of the first $n$ customers, where $D_{j:n}$ denotes the set of the
$n$ customers seated at a table with common rank $j$.

\Enumerate \item[(i)]Given this configuration, the next customer
$n+1$ is seated at an occupied table $D_{j:n}$, denoting that
customer $n+1$ is equivalent to the $j$th largest seated
customers, with probability, \Eq \frac{p_{j:n}}{l(n)}\int_{\Xcr}
K(Y_{n+1}|x)\pi(dx|D_{(j:n)}) \label{seatjold} \EndEq for
$j=1,\ldots, n(\p)$. \item[(ii)]Otherwise, the probability that
customer $n+1$ is new and is the $j$th largest among $n(\p)+1$
possible ranks is, \Eq \frac{q_{j:n}}{l(n)}\int_{\Xcr}
K(Y_{n+1}|x)H(dx) \label{seatjnew} \EndEq for $j=1,\ldots,
n(\p)+1$. \EndEnumerate Similar to the gWCR SIS algorithms~[see
Ishwaran and James~(2003, Lemma 2)], by appealing to the product
rule of probability, repeating this procedure for customers
$\{1,\ldots,n\}$, produces a draw of $\m$ from a density of $\m$
depending on $\Y$, say $q(\m)$, that satisfies the relationship $$
L(\m)q(\m)=
\pi_{\rho}(\m)\prod_{j=1}^{n(\p)}\int_{\Xcr}\prod_{i\in
D_{j}}k(Y_{i}|x)H(dx)$$ where $L(\m)=\prod_{i=1}^{n}l(i-1)$. Hence
for any functional, $h(\m)$ it follows that \Eq
\sum_{\p}\sum_{\m\in S_{n(\p)}}h(\m)\pi_{\rho}(\m|\Y)=
\frac{\sum_{\p}\sum_{\m \in S_{n(\p)}}h(\m)L(\m)q(\m)}
{\sum_{\p}\sum_{\m \in S_{n(\p)}}L(\m)q(\m)}. \label{realpost}
\EndEq If the functional $h(\m)$ has a closed form, such as the
predictive density $\E[f(y|P)|\m,\Y]=f(y|\m,\Y)$, then one
approximates~\mref{realpost} by using the rules in~\mref{seatjold}
and ~\mref{seatjnew} to draw $\m$. Repeating this procedure say
$B$ times, results in iid realizations say $(\m_{(b)})$ for
$b=1,\ldots, B$ and one can approximate~\mref{realpost} by
$$
\frac{\sum_{b=1}^{B}h(\m_{(b)})L(\m_{(b)})}
{\sum_{b=1}^{B}L(\m_{(b)})}.
$$
When the kernels $K$ are set to $1,$  this procedure reduces to
that described in James~(2003, 2006) producing an exact draw from
$\pi_{\rho}(\m).$ For more intricate models one can incorporate a
draw from the unique values $X^{*}_{1},\ldots, X^{*}_{n(\p)}$
which has the same distribution that arises for the Dirichlet
process. One can also incorporate draws from the posterior
distribution of $P_{\rho,H}(dx)$ which is described in
James~(2006). Otherwise it is a simple matter to modify all the
computational procedures discussed in Ishwaran and James (2003,
section 4).
\subsection{Normal Mixture example}
One of the most studied and utilized Bayesian mixture models is
the Normal mixture model, specified by the choice of \Eq
f_{\sigma}(y|P)=\int_{-\infty}^{\infty}\phi_{\sigma}(y-x)P(dx)
\label{Nmix}\EndEq where
$$ \phi_{\sigma}(z)=\frac{1}{\sqrt{2\pi}\sigma}\exp\(
-\frac{{z}^{2}}{2\sigma^{2}}\) $$ is a Normal density, which is a
natural candidate for density estimation. In the case of the
Dirichlet process, this model was introduced by Lo~(1984) and
popularized by the development of feasible computational
algorithms in Escobar~(1994) and Escobar and West~(1995). Suppose
that $(Y_{i})$ are iid with true density $f_{0}$, a recent result
of Lijoi, Pr\"unster and Walker~(2005) shows that
$f_{\sigma}(\cdot|P)$ in~\mref{Nmix} based on very general random
probability measures, and a suitable prior distribution for
$\sigma$, have posterior distributions that are strongly
consistent in terms of estimating the unknown density $f_{0}$
under rather mild conditions. In particular their result validates
the use of rather arbitrary NTR species sampling models in this
context with the classical choice of $H$ set to be a Normal
distribution with mean 0 and variance $A.$  Here setting
$\sigma=\sqrt{\theta}$ one has
$$
K(Y_{i}|X_{i})=\frac{1}{\sqrt{2\pi
\theta}}\exp\(\frac{1}{2\theta}{(Y_{i}-X_{i})}^{2}\).
$$
Using these specifications we present the details of the proposed
algorithm.
\begin{enumerate}
\item[(i)] Customer $n+1$ is seated to a new table and assigned rank
$j$ among $n(\p)+1$ possible ranks with probability
$$
\frac{q_{j:n}}{\lambda_{\theta}(n+1)}\frac{1}{\sqrt{2\pi(\theta+A)}}\exp\(-\frac{Y^{2}_{n+1}}{2(\theta+A)}\)
$$
\item[(ii)] Customer $n+1$ is seated to an existing table and is assigned
rank $j$ with probability
$$
\frac{p_{j:n}}{\lambda_{\theta}(n+1)}\sqrt{\frac{\theta+Ad_j}{2\pi\theta[\theta+A(d_{j}+1)]}}
\exp\[-\frac{1}{2\theta}\(Y^{2}_{n+1}-\frac{A\sum_{i\in
D_{j}}Y_{i}+Y_{n+1}}{\theta+A(d_{j}+1)}+\frac{A\sum_{i\in
D_{j}}Y_{i}}{\theta+Ad_{j}}\)\]
$$
\item[(iii)] Additionally each $X^{*}_{j}|\Y,\m,\theta$ is normally
distributed with parameters
$$
\frac{1}{\sigma_{j}}=\frac{d_{j}}{\theta}+\frac{1}{A}{\mbox { and
}} \mu_{j}=\frac{\sigma_{j}}{\theta}\sum_{i\in D_{j}}Y_{i}.
$$
\end{enumerate}
$\lambda_{\theta}(n+1)$ is the appropriate normalizing constant
which is a special case of $l(n).$\Remark For comparison, the
setup and notation we use is similar to that used in Ishwaran and
James~(2003, 6.1) which is based on weighted Chinese restaurant
sampling of partitions $\p$. \EndRemark
 \Section{Concluding Remarks} We have given a brief
account of how one can use Kjell Doksum's NTR models to create a
new class of species sampling random probability measures which
can be applied to complex mixture models. These models exhibit
many features of the NTR models, in terms of clustering behavior,
but as we have shown are simpler to use. Ideally one would like to
describe parallel schemes for the more complex Spatial NTR models.
However, this constitutes a considerably more involved study which
we shall report elsewhere. More details can be found in
James~(2003, 2006) where explicit examples can be easily
constructed.

The representation in~\mref{expref} is important as it connects
NTR processes to a large body of work on exponential functionals
of L\'evy processes which have applications in many fields
including physics and finance. For a recent survey see Bertoin and
Yor~(2005). Some recent works which exploit this representation
and are directly linked to NTR processes are Epifani, Lijoi and
Pr\"unster~(2003) and
    James~(2003, 2006). Additionally, outside of a Bayesian context, there is a notable body of recent work
    which has some overlaps with James~(2003, 2006) and hence NTR
    processes by Gnedin and Pitman~(2005a) and subsequent papers Gnedin and Pitman~(2005b),
     Gnedin and Pitman and Yor~(2005) and Gnedin, Pitman and Yor~(2006).
    Although outside of a specific Bayesian context these papers
    contain results which are relevant to statistical analysis
    such as results related to the behavior of the number of ties
    $n(\p).$ The fact that these models arise from
    different considerations and different points of emphasis
    attests to their rich nature. We are quite interested to see
    what future connections will be made.
\vskip0.2in

\centerline{\Heading References} \vskip0.4in \tenrm
\def\smc{\tensmc}
\def\sl{\tensl}
\def\bf{\tenbold}
\baselineskip0.15in

\Ref \by    Antoniak, C. E. \yr    1974 \paper Mixtures of
Dirichlet processes with applications to Bayesian nonparametric
problems \jour  \AnnStat \vol   2 \pages 1152-1174 \EndRef

\Ref \by    Bertoin, J. and Yor, M. \yr    2005 \paper Exponential
functionals of L\'evy processes\jour  Probab. Surv. \vol 2 \pages
191-212\EndRef

\Ref \by    Blackwell, D. and MacQueen, J. B. \yr 1973 \paper
Ferguson distributions via P\'olya urn schemes \jour \AnnStat \vol
1 \pages 353-355 \EndRef

\Ref \by    Doksum, K. A. \yr    1974 \paper Tailfree and neutral
random probabilities and their posterior distributions \jour
\AnnProb \vol   2 \pages 183-201 \EndRef

\Ref \yr    2004 \by    Doksum,~K.~A. and James,~L.~F. \paper On
spatial neutral to the right processes and their posterior
distributions. In Mathematical Reliability: An Expository
Perspective, Editors: Mazzuchi, Singpurwalla and Soyer.
International Series in Operations Research and Management
Science. Kluwer Academic Publishers \EndRef

\Ref \by    Epifani, I., Lijoi, A., and Pruenster, I. \yr    2003
\paper Exponential functionals and  means of neutral to the right
priors \jour Biometrika \vol 90 \pages 791-808\EndRef

\Ref \by Escobar, M.D. \yr 1994 \paper Estimating normal means
with the Dirichlet process prior \jour \JASA \vol 89 \pages
268-277 \EndRef

\Ref \by Escobar, M.D. and West, M. \yr 1995 Bayesian density
estimation and inference using mixtures. \jour JASA \vol 90 \pages
 577-588 \EndRef

\Ref \by    Ewens, W. J. \yr    1972 \paper The sampling theory of
selectively neutral alleles \jour  Theor. Popul. Biol. \vol   3
\pages 87-112 \EndRef

\Ref \by    Ferguson, T. S. \yr    1973 \paper A Bayesian analysis
of some nonparametric problems \jour  \AnnStat \vol   1 \pages
209-230 \EndRef

\Ref \by     Ferguson, T. S. \yr     1974 \paper  Prior
distributions on spaces of probability measures \jour   \AnnStat
\vol    2 \pages  615-629 \EndRef

\Ref \by     Ferguson, T. S. and Phadia, E. \yr     1979 \paper
Bayesian nonparametric estimation based on censored data \jour
\AnnStat \vol    7 \pages  163-186 \EndRef

\Ref \by     Freedman, D. A. \yr     1963 \paper  On the
asymptotic behavior of
  Bayes estimates in the discrete case
\jour   Ann. Math. Statist. \vol    34 \pages  1386-1403 \EndRef

\Ref \by Gnedin, A. V. \yr 2004 \paper Three sampling
formulas\jour Combin. Probab. Comput. \vol 13 \pages
185-193\EndRef \Ref \by Gnedin, A. V. and Pitman, J. \yr  2005a
\paper Regenerative composition structures \jour \AnnProb \vol 33
\pages 445-479\EndRef

\Ref \by      Gnedin, A. V. and Pitman, J. \yr  2005b  \paper
Self-similar and Markov composition structures. In Representation
Theory, Dynamical Systems, Combinatorial and Algorithmic Methods.
Part 13, A. A. Lodkin editor. Zapiski Nauchnyh Seminarov POMI,
Vol. 326, PDMI, 59-84\EndRef

\Ref \by      Gnedin, A. V. and Pitman, J. and Yor, M. \yr  2005
\paper Asymptotic laws for regenerative compositions: gamma
subordinators and the like \jour Probab. Th. and Rel. Fields.
Published online November 2005\EndRef

\Ref \by      Gnedin, A. V. and Pitman, J. and Yor, M. \yr  2006
\paper Asymptotic laws for compositions derived from transformed
subordinators \jour \AnnProb \vol 34\EndRef

\Ref \by      Hjort, N. L. \yr      1990 \paper   Nonparametric
Bayes estimators based on Beta processes in
         models for life history data
\jour    \AnnStat \vol         18 \pages    1259-1294 \EndRef

\Ref \yr    2001 \by    Ishwaran,~H. and James,~L.~F. \paper Gibbs
sampling methods for stick-breaking priors \jour Journal of the
American Statistical Association \vol 96 \pages 161-173 \EndRef

\Ref \yr    2003 \by    Ishwaran,~H. and James,~L.~F. \paper
Generalized weighted Chinese restaurant processes for species
       sampling mixture models
\jour  Statistica Sinica \vol   13 \pages 1211-1235 \EndRef

\Ref \yr    2003 \by    James,~L.~F. \paper Poisson calculus for
spatial neutral to the right processes(Big version).
\\arXiv:math.PR/0305053. Available at
http://arxiv.org/abs/math.PR/0305053 \EndRef

\Ref \yr    2006 \by    James,~L.~F. \paper Poisson calculus for
spatial neutral to the right processes\jour \AnnStat \vol
34\EndRef

\Ref \by    Kim, Y. \yr    1999 \paper  Nonparametric Bayesian
estimators for counting processes \jour   \AnnStat \vol    27
\pages  562-588 \EndRef

\Ref \by Lijoi, A., Pr\"unster, I. and Walker, S.G. \yr 2005
\paper On consistency of nonparametric normal mixtures for
Bayesian density estimation\jour \JASA \vol 100 \pages
1292-1296\EndRef

\Ref \by Lo, A. Y. \yr 1993 \paper A Bayesian bootstrap for
censored data \jour  \AnnStat \vol 21 \pages 100-123 \EndRef

\Ref \by    Lo, A. Y. \yr    1984 \paper On a class of Bayesian
nonparametric estimates: I.  Density
       Estimates
\jour  \AnnStat \vol   12 \pages 351-357 \EndRef

\Ref \by Lo, A.Y., Brunner, L.J. and Chan, A.T. \yr 1996 \paper
Weighted Chinese restaurant processes and Bayesian mixture
       model. Research Report Hong Kong University of
       Science and Technology
\EndRef

\Ref \by M\"uller, P, and  Quintana, F. A.\yr 2004 \paper
Nonparametric Bayesian data analysis \jour Statist. Sci. \vol 19
\pages 95-110 \EndRef

\Ref \by Pitman, J. \yr    1996 \paper Some developments of the
Blackwell-MacQueen urn scheme. In Statistics, Probability and Game
Theory T.S. Ferguson, L.S. Shapley and J.B. Macqueen editors, IMS
Lecture Notes-Monograph series, Vol 30, pages 245-267 \EndRef

\Ref \by    Pitman, J. and Yor, M. \yr    1997 \paper The
two-parameter Poisson-Dirichlet distribution derived from
       a stable subordinator
\jour  \AnnProb \vol   25 \pages 855-900 \EndRef

\Ref \by    Walker, S. and Muliere, P. \yr    1997 \paper
Beta-Stacy processes and a generalization of the P{\'o}lya-urn
scheme \jour   \AnnStat \vol    25 \pages  1762-1780 \EndRef


\smc

\Tabular{ll}

Lancelot F. James\\
The Hong Kong University of Science and Technology\\
Department of Information Systems and Management\\
Clear Water Bay, Kowloon\\
Hong Kong\\
\rm lancelot\at ust.hk\\

\EndTabular

\end{document}